\documentclass[leqno]{amsart}
\usepackage[T1]{fontenc}
\RequirePackage{amsthm,amsmath,amsfonts}
\RequirePackage[colorlinks,citecolor=blue,urlcolor=blue]{hyperref}
\newtheorem{theorem}{Theorem}[section]
\newtheorem{lemma}[theorem]{Lemma}
\newtheorem{proposition}[theorem]{Proposition}

\newtheorem{corollary}[theorem]{Corollary}
\theoremstyle{definition}
\newtheorem{definition}[theorem]{Definition}

\theoremstyle{remark}
\newtheorem{remark}[theorem]{Remark}
\numberwithin{equation}{section}

\begin{document}
\title{Vol-of-vol expansion for (rough) stochastic volatility models}
\author{Ozan Akdo\u{g}an}
\begin{abstract}
We introduce an asymptotic small noise expansion, a so called \emph{vol-of-vol expansion}, for potentially infinite dimensional and rough stochastic volatility models. Thereby we extend the scope of existing results for finite dimensional models and validate claims for infinite dimensional models. Furthermore we provide new, explicit (in the sense of non-recursive) representations of the so-called push-down Malliavin weights that utilizes a precise understanding of the terms of this expansion.
\end{abstract}
\date{\today}
\maketitle

\section{Introduction}
With the transition from the Black \& Scholes model to stochastic volatility models the lack of analytic pricing formulae made the calibration against the implied volatility surface a computational challenge. Asymptotic expansions (see \cite{friz2015large} for an overview) provide fast methods for approximate pricing and hence are popular for calibration. As there is now strong evidence from option prices (cf. \cite{alos2007short} and \cite{fukasawa2011asymptotic}) as well as from historical prices (cf. \cite{gatheral2018volatility}) that stochastic \emph{volatility is rough}, the benefits of asymptotic expansions are even more pronounced.

The so-called \emph{vol-of-vol expansion} is particularly popular, as the leading term agrees with the Black \& Scholes price and higher order terms are usually given in closed form. An early version of this expansion for one dimensional stochastic volatility models is given in \cite{lewis2000option} and a generalization to finite-dimensional stochastic volatility models in \cite{TT2012}. A further formal generalization came with the popular Bergomi-Guyon expansion \cite{bergomi2012stochastic} (henceforth \emph{BG expansion}) which is formulated in terms of potentially infinite dimensional forward variance models and was incremental in the introduction of the rough Bergomi model (cf. \cite{bayer2016pricing}). Although no conditions for the validity of this expansion are provided, it is easy to see that the terms agree with the ones in \cite{TT2012} for finite-dimensional stochastic volatility models. In \cite{alos2018exponentiation} a series representation, the so-called \emph{forest expansion}, is introduced and shown that, when appropriately truncated, the terms agree with those of the BG expansion for the rough Heston model but not for the rough Bergomi model. However, there is no convergence rate provided. Hence precise conditions for the validity of the BG expansion for the infinite dimensional (and potentially rough) case remain to the best of our knowledge unknown. 

The main contributions of this paper will be precise conditions for this expansion to hold and new explicit representations that improve on the recursive representations given in \cite{TT2012}. We will also show the relation to the forest expansion given in \cite{alos2018exponentiation} and prove (under some mild conditions) that a necessary condition for this expansion to agree with the BG expansion is that the underlying forward variance model is affine.

The rest of this paper is organized as follows. In Section \ref{secMain} we introduce the setting and state the main result of this paper, namely conditions for the vol-of-vol expansion to hold and new explicit representations of the terms of the expansion. After providing some examples of forward variance models to which this expansion applies, we will compare it in more detail to the vol-of-vol expansions mentioned above. In the remaining sections we will prove the main result by utilizing the Watanabe-Yoshida theory, i.e. the Yoshida extension \cite[Theorem 2.2]{YN1992} of Watanabe's Theorem \cite[Theorem 2.3]{WS1987}. This Yoshida extension will be crucial as in the original version of this theorem the non-degeneracy of the Malliavin covariance matrix of the forward variance is necessary. While for finite dimensional stochastic volatility models the H\"{o}rmander condition can be utilized, in the infinite dimensional case this is not that straight forward (cf. \cite{gerasimovics2018h} for semi linear SPDEs with negative definite self-adjoint operators). With the Yoshida extension it will be sufficient to have a local non-degeneracy condition. As the setting \cite[Theorem 2.2]{YN1992} is tailored for finite dimensional SDEs, we prove in Section \ref{sWTA} a version of this theorem that is suitable for our infinite dimensional SPDE setting. The results in this section are formulated in a quite general way and could be used for example to tackle the expansion of the term structure of forward interest rates as considered in \cite[Section 4]{KT2003} directly in the infinite dimensional setting, i.e. without an initial discretization to arrive at the finite dimensional setting. In Section \ref{SNE} we will then provide conditions such that the truncated version of the weak Taylor expansion validates the vol-of-vol expansion for a large class of forward variance models, which include infinite dimensional but not yet rough forward variance models. In Section \ref{PDMW} we will give explicit representations of the terms appearing in the expansion. Based on this representations and a change of limits argument we will show in Section \ref{RFVM} that this expansion holds true (essentially) whenever it exists, which in particular is the case for rough models.
\section{Main result}\label{secMain}
\subsection{Setting and notation}\label{secnot}
For fixed $T>0$, we assume as given a filtered probability space $(\Omega,\mathcal{F},(\mathcal{F}_t)_{t\in [0,T]},\mathbb{P})$ accommodating a $2$-dimensional Brownian motion $\beta=(\beta^1,\beta^2)$, with $\mathbb{P}$ being a local martingale measure for the market of stock prices and variance swaps. For a given Hilbert space $\mathcal{H}$ we will denote by $C_{b}^\infty(\mathcal{H})$ the set of smooth functions on $\mathcal{H}$ with bounded derivatives and by $C_{bb}^\infty(\mathcal{H})$ the set of functions in $C_b^\infty(\mathcal{H})$ that are bounded. Further, we denote by $C_p^\infty(\mathcal{H})$ the set of smooth functions on $\mathcal{H}$ such that the functions and its derivatives have polynomial growth. If $\mathcal{H}=\mathbb{R}$, we will write $C_b^\infty$, $C_{bb}^\infty$ and $C_p^\infty$. 

We will mainly use the \emph{term structure} Hilbert space $H$ that was introduced in \cite[Chapter 5]{DF2001}. Elements of $H$ are bounded and absolutely continuous and the point evaluation is a continuous functional. The multiplication operator $u\mapsto m(u,g)$ is a linear continuous operator on $H$ if and only if $g\in H$ (see \cite[Lemma II.5.6]{OA2016}). 

Let $\sigma$ be Lipschitz continuous on $H$ and satisfy a linear growth condition, then for $g\in H$ the volatility of forward variance vector field on $H$, given by $u\mapsto m(\sigma(u),g)$ is also Lipschitz continuous and satisfies a linear growth condition. For $\epsilon\in [-1,1]$ we denote the $H$ valued forward variance process in Musiela's parameterization (see \cite{musiela1993stochastic}) by $u_t^\epsilon$ and assume that it is given by the continuous mild solution of 
\begin{equation*}\label{ForwV}
du_t^\epsilon=\frac{d}{dx}u_t^\epsilon\,dt+\epsilon\, m(\sigma(u^\epsilon_t),g)\,d\beta^1_t,\quad u^\epsilon_0=u,
\end{equation*}
i.e., it satisfies  (see \cite[Chapter 7]{da2014stochastic})
\begin{equation}\label{ForwVmild}
u_t^\epsilon=S_tu +\epsilon\,\int_0^t S_{t-s}m(\sigma(u^\epsilon_s),g)\,d\beta^1_s.
\end{equation}
Here $S_t$ denotes the shift-operator acting on real-valued functions $f$ by $S_tf(x) = f(t+x)$. The corresponding log-price $X_t^\epsilon$ is then given by
\begin{equation}\label{logPrice}
dX_t^\epsilon=-\frac{1}{2}u_t^\epsilon(0)\,dt+\sqrt{u_t^\epsilon(0)}\,d\beta_t,\quad X^\epsilon_0=x,
\end{equation}
where $d\beta_t:=(\rho\, d\beta^1_t+\sqrt{1-\rho^2}\, d\beta_t^2)$ with correlation $\rho\in [-1,1]$. 

For $\Sigma>0$ we will recall (cf. \cite[Exercise 1.1.1]{nd2006}) that the generalized Hermite polynomials $H_n(x,\Sigma)$, $n\geq 0$, satisfy 
\begin{equation*}\label{Hn1}
H_n(x,\Sigma)=\frac{(-\Sigma)^n}{n!}e^{\frac{x^2}{2\Sigma}}\frac{d^n}{dx^n}e^{-\frac{x^2}{2\Sigma}},\quad \frac{d}{dx}H_n(x,\Sigma)=H_{n-1}(x,\Sigma),\quad n\geq 1.
\end{equation*}
We will write $H_k(T):=H_k(Y_T,\Sigma_T)$ for $k\geq 0$, where
\begin{equation*}\label{YSigma}
Y_T:=\int_0^T \sqrt{u(t_1)}\,d\beta_{t_1}\quad\text{and}\quad \Sigma_T:=\int_0^Tu(t_1)\,dt_1,
\end{equation*}
denote the martingale part and quadratic variation of $X_T^0$. Finally for later reference we mention the representation 
\begin{equation*}\label{Hito}
H_{m}(T)=\int_0^TH_{m-1}(t_m)\sqrt{u(t_m)}\, d\beta_{t_m}\quad \text{for}\quad m\geq 1,
\end{equation*}
that follows from \cite[Exercise 1.1.1]{nd2006} and \cite[(1.27)]{nd2006}. We will further denote 
\begin{equation*}\label{ui}
U^{(i)}:=\frac{\partial^i}{\partial \epsilon^i}\big\vert_{\epsilon=0}U^\epsilon\quad \text{ and }\quad \sqrt{U}^{(i)}:=\frac{\partial^i}{\partial \epsilon^i}\big\vert_{\epsilon=0}\sqrt{U^\epsilon},
\end{equation*}
and by $B_{n,k}(T):=B_{n,k}(X_T^{(1)},...,X_T^{(n-k+1)})$ the partial Bell polynomial 
\begin{eqnarray}\label{partBell}
B_{n,k}(T)=\sum_{\mathbf{j}\in T(n,k)}\frac{n!}{\prod_{i=1}^{n-k+1}j_i!}\prod_{i=1}^{n-k+1}\Big(\frac{X_T^{(i)}}{i!}\Big)^{j_i},
\end{eqnarray}
with $T(n,k)$ being the set of tuples $(j_1,...,j_{n-k+1})$ of non-negative integers satisfying
\begin{equation*}\label{tnk}
\sum_{i=1}^{n-k+1}j_i=k\quad \text{ and }\quad \sum_{i=1}^{n-k+1}ij_i=n.
\end{equation*} 
In the following we set $B_{0,0}(T)=1$ and $B_{k,0}(T)=B_{0,k}(T)=0$ if $k\geq 1$.
\subsection{Main result}
We notice that for $g\in L^2([0,T])$ the spot-variance $u_t^\epsilon(0)$ and the log-price \eqref{logPrice} are still well-defined.
\begin{theorem}\label{MainThm}Let $u\in H$ be a strictly positive initial curve, $\sigma\in C^\infty_b(H)$ and $f$ being in $C^\infty_p$ or bounded. If either $g\in H$, or $g\in L^2([0,T])$ and $f$ being also Lipschitz continuous, then the following weak Taylor expansion holds up to arbitrary order $p\geq 2$,
\begin{equation*}\label{MEXP}
\mathbb{E}[f(X_T^\epsilon)]=\sum_{i=0}^p\sum_{k=0}^i\sum_{l=k}^{i+2k}\frac{\epsilon^i}{i!}\mathbb{E}[B_{i,k}(T)H_{l-k}(T)]\frac{d^l}{dx^l}\mathbb{E}[f(X_T^0)]+o(\epsilon^p)\,\,\,\text{as }\epsilon\rightarrow 0.
\end{equation*}
\begin{proof}
This follows from Propositions \ref{propSNE}, \ref{expREep} and \ref{roughWTX}.
\end{proof}
\end{theorem}
Thus the terms of the expansion at $p$-th order are given by expected values of polynomials in $X^{0},X^{(1)},...,X^{(p)}$ as well as derivatives up to $3p$-th order of the Black \& Scholes price $\mathbb{E}[f(X_T^0)]$ with respect to the initial log-price.
\subsection{Examples}\label{exeps}
Our main examples that admit the vol-of-vol expansion are the generalized stochastic volatility models with affine drift and the generalized Bergomi model (all details can be found in \cite{OA2016}). The former class of models correspond to \eqref{ForwVmild} with $\sigma(u)=\tilde{\sigma}(u(0))1$ with $1\in H$ and the latter to $\sigma(u)=u$ while in both cases $g\in L^2([0,T])$. Popular examples for the former class of models are $\tilde{\sigma}(x)=\sqrt{x}$ (see Remark \ref{sqremark}) and $\tilde{\sigma}(x)=x$ for the Heston and the GARCH model, respectively. 

Both model classes admit a finite dimensional realization if (and essentially only if) $g(x)=\phi e^{-bx}$ for positive numbers $\phi$ and $b$ (and hence $g\in H$). In this case the former class of models correspond to (Hull-White extended) stochastic volatility models in which $v_t:=u_t(0)$ satisfies
\begin{equation*}
dv_t=b(\theta-v_t)\,dt+\phi \tilde{\sigma}(v_t)\,d\beta_t^1,\quad v_0=v,
\end{equation*}
and the generalized Bergomi model corresponds in this case to the conventional Bergomi model. The rough versions of this models are given for the choice $g(x)=\phi x^{-\gamma}$ with $\phi>0$ and Hurst parameter $H=\frac{1}{2}-\gamma \in (0,\frac{1}{2})$ (and hence $g\in L^2([0,T])$). As we will show now, this models have the pleasant property that they resemble the term structure of the at-the-money volatility skew very well. For this we recall that the corresponding expansion of the implied volatility (cf. \cite[Theorem 3.2]{TT2012} or \cite[Section 2]{bergomi2012stochastic}) at first order, is given by 
\begin{equation*}
\sigma^\epsilon=\sigma^{(0)}+\epsilon\sigma^{(1)}+o(\epsilon),\quad \text{as } \epsilon\rightarrow 0,
\end{equation*}
with $\sigma^{(0)}=\sqrt{\frac{1}{T}\int_0^Tu(t_1)\,dt_1}$ and
\begin{equation*}\label{sigma1}
\sigma^{(1)}(k,T)=\Big(\frac{1}{2}-\frac{k}{\int_0^Tu(t_1)\,dt_1}\Big)\frac{\rho\int_0^T\int_0^{t_1} m(\sigma(u_{t_2}^0),g)(t_1-t_2)\sqrt{u(t_2)}\,dt_2\,dt_1}{2\sqrt{T\int_0^Tu(t_1)\,dt_1}},
\end{equation*}
where $x$ denotes the initial log-price, $k=\log(e^x/K)$ the log-strike of the call option expiring in $T$. Then the term structure of the at-the-money volatility skew satisfies at first order (cf. \cite{bayer2016pricing}) for a flat initial curve $u=1$  is given by
\begin{equation*}
\psi(T;g):=\Big\vert \frac{\rho\sigma(1)\int_0^T\int_0^{t_1} g(t_1-t_2)\,dt_2\,dt_1}{2T^2}\Big\vert\propto\Big\vert \frac{\int_0^T\int_0^{t_1} g(t_1-t_2)\,dt_2\,dt_1}{2T^2}\Big\vert\label{atmskew}.
\end{equation*}
Hence in the rough version of the models we have $\psi(T;g)\propto T^{-\gamma}$ which is just the empirically observed form if $\gamma\approx 0.4$ and hence $H\approx 0.1$. 
\subsection{Related literature}
The expansion introduced above is a generalization of the results in \cite[Section 3.2]{TT2012} to infinite dimensional forward variance models. Moreover, while coming to the same results, the Representation \eqref{MEXP} is explicit as opposed to the representations given in \cite[Theorem 3.1.]{TT2012} which are given in terms of recursively defined Malliavin weights.

Besides the different parameterizations of the forward variance (\emph{time of maturity} vs \emph{time to maturity}) our expansion agrees with the BG expansion. However, while in \cite{bergomi2012stochastic} the asymptotic property is only assumed we provide precise conditions. Further, our explicit representation allows us to sharpen the interpretation of the terms $C_0^{xu}(u), C_0^{uu}(u)$ and $C_0^\mu(u)$ which are introduced in \cite{bergomi2012stochastic}. 

With $F^l:=\frac{d^l}{dx^l}\mathbb{E}[f(X_T^0)]$, the expansion up to second order is given by 
\begin{equation} \label{MEXP2}
\begin{cases}
\mathbb{E}[f(X_T^\epsilon)]&=\mathbb{E}[f(X_T^0)]+\epsilon\,\sum_{l=2}^3\mathbb{E}[X_T^{(1)}H_{l-1}(T)]F^l\\
&+\frac{\epsilon^2}{2}\sum_{l=2}^4(\mathbb{E}[X_T^{(2)}H_{l-1}(T)+(X_T^{(1)})^2H_{l-2}(T)])F^l\\
&+\frac{\epsilon^2}{2}\sum_{l=5}^6\mathbb{E}[(X_T^{(1)})^2H_{l-2}(T)]F^l+o(\epsilon^2)\quad \text{as } \epsilon\rightarrow 0,
\end{cases}
\end{equation}
where we used that $\mathbb{E}[X_T^{(1)}H_{0}(T)]=\mathbb{E}[X_T^{(2)}H_{0}(T)]=0$. In the following we write $\Sigma(u):=m(\sigma(u),g)$. Then from \eqref{xmjm} and \eqref{um} it is straight-forward to compute the individual terms of \eqref{MEXP2}. In the notation of \cite{bergomi2012stochastic} we have
\begin{equation*}
\begin{cases}\label{cxu}
&C_0^{xu}(u)=\rho\int_0^T\int_0^{t_1}\Sigma(u_{t_2}^0)(t_1-t_2)\sqrt{u(t_2)}\,dt_2\,dt_1\\
&C_0^{uu}(u)=2\int_0^T\int_0^{t_1} \int_0^{t_2}\Sigma(u_{t_3}^0)(t_1-t_3)\Sigma(u_{t_3}^0)(t_2-t_3)\,dt_3\,dt_2\,dt_1\\
&C_0^\mu(u)=\rho^2\Big(\int_0^T \int_0^{t_1}\int_0^{t_2}d\Sigma(u_{t_2}^0)\Sigma(u_{t_3}^0)(t_1-t_3)\sqrt{u(t_3)}\sqrt{u(t_2)}\,dt_3\,dt_2 \,dt_1\\
&\quad +\frac{1}{2}\int_0^T \int_0^{t_1}\int_0^{t_2}\frac{\Sigma(u_{t_2}^0)(t_1-t_2)}{\sqrt{u(t_2)}} \Sigma(u^0_{t_3})(t_2-t_3)\sqrt{u(t_3)}\,dt_3 \,dt_2\,dt_1\Big),
\end{cases}
\end{equation*}
and the terms of \eqref{MEXP2} at first order are given by 
\begin{eqnarray*}
&&\mathbb{E}[X_T^{(1)}H_1(T)]F^2=-\frac{1}{2}C_0^{xu}(u)F^2,\quad \mathbb{E}[X_T^{(1)}H_2(T)]F^3=\frac{1}{2}C_0^{xu}(u)F^3,
\end{eqnarray*}
and at second order by
\begin{eqnarray*}
&&(\mathbb{E}[X_T^{(2)}H_1(T)+(X_T^{(1)})^2H_0(T)])F^2=\frac{1}{4}C_0^{uu}(u)F^2,\\
&&(\mathbb{E}[X_T^{(2)}H_2(T)+(X_T^{(1)})^2H_1(T)]) F^3=(-C_0^\mu(u)-\frac{1}{2}C_0^{uu}(u))F^3,\\
&&(\mathbb{E}[X_T^{(2)}H_3(T)+(X_T^{(1)})^2H_2(T)])F^4\\
&&\quad\quad\quad\quad\quad\quad\quad\quad\quad=(C_0^\mu(u)+\frac{1}{4}C_0^{uu}(u)+\frac{1}{4}(C_0^{xu}(u))^2)F^4,
\end{eqnarray*}
as well as
\begin{eqnarray*}
&&\mathbb{E}[(X_T^{(1)})^2H_3(T)]F^5=-\frac{1}{8}(C_0^{xu}(u))^2F^5,\\
&&\mathbb{E}[(X_T^{(1)})^2H_4(T)]F^6=\frac{1}{4}(C_0^{xu}(u))^2F^6.
\end{eqnarray*}
In \cite{bergomi2012stochastic} the terms $C_0^{xu}(u)$ and $C_0^{uu}(u)$ are introduced as \emph{the integrated spot / variance covariance function} and \emph{the integrated variance / variance covariance function}, respectively, and formally defined as
\begin{eqnarray*}
&&C_0^{x\xi}(\xi)=\int_0^T\int_{s}^T\frac{\mathbb{E}[\frac{dS_s}{S_s}d\xi_s^x]}{ds}\,dx\,ds\\
&&C_0^{\xi\xi}(\xi)=\int_0^T\int_s^T\int_s^T\frac{\mathbb{E}[d\xi_s^xd\xi_s^y]}{ds}\,dx\,dy\,ds,
\end{eqnarray*}
where $\xi$ is the initial curve and $\xi_s^x$ denotes the forward variance at time $s$ for time $x$ and $S$ denotes the stock price process. It is related to the forward variance $u$ in Musiela's parameterization by $\xi_s^x=u_s(x-s)$. In our notation these terms are 
\begin{eqnarray*}
C_0^{xu}(u)&=&\int_0^T\mathbb{E}[\langle u^{(1)}(0),X^0\rangle_{t_1}]\,dt_1,\\
C_0^{uu}(u)&=&4\int_0^T\int_0^{t_1}\mathbb{E}[\langle u^{(1)}(t_1-t_2), u^{(1)}(0)\rangle_{t_2} ]\,dt_1\,dt_2,
\end{eqnarray*}
and thus showing that the first term is the integrated covariation between $u^{(1)}(0)$ and $X^{(0)}$ and the second term the integrated auto-covariation of $u^{(1)}$. 

Finally we mention the expansions given in \cite{alos2018exponentiation} and \cite{gulisashvili2019higher} which are generalizations of \cite{alos2012decomposition}. In the latter a decomposition formula for option prices in the Heston model is introduced and based on that approximations proposed that are of first order in vol-of-vol. In \cite{gulisashvili2019higher} this approximation for the Heston model is extended to third order in vol-of-vol. In \cite{alos2018exponentiation} the decomposition is extended to a series representation and according to \cite[Section 4]{alos2018exponentiation} the terms of the BG expansion $C_0^{x,u}(u)$, $C_0^{u,u}(u)$ and $C_0^{\mu}(u)$ correspond to the terms (of their \emph{forest expansion}) $(X\diamond M)$, $(M\diamond M)$ and $(X\diamond (X\diamond M))$, respectively, while noting that equality does not always hold. In the current notation, $X_t=X_t^\epsilon$ and 
\begin{equation*}
M_t=M_t^\epsilon=\epsilon\int_0^t\int_{t_1}^T m(\sigma(u_{t_1}^\epsilon),g)(t_2-t_1)\,dt_2\,d\beta^1_{t_1}.
\end{equation*}
Then the first term satisfies
\begin{eqnarray*}
&&(X\diamond M)_0(T)=\epsilon\,\rho\,\int_0^T\int_{t_1}^T\mathbb{E}[m(\sigma(u_{t_1}^\epsilon),g)(t_2-t_1)\sqrt{u_{t_1}^\epsilon(0)}]\,dt_2\,dt_1
\end{eqnarray*}
and hence when comparing it (after an application of the Fubini Theorem and setting $\Sigma(u)=m(\sigma(u),g)$) with the corresponding term in \eqref{cxu} we see 
\begin{eqnarray*}
&&(X\diamond M)_0(T)=\epsilon\,\mathbb{E}[C_0^{x,u}(u^\epsilon)]
\end{eqnarray*}
and similarly $(M\diamond M)_0(T)=\epsilon^2\,\mathbb{E}[C_0^{u,u}(u^\epsilon)]$ and so on for the higher order terms. Surprisingly, for the affine forward variance models introduced in \cite{gatheral2018affine} (corresponding to the stochastic volatility models with affine drift given in \ref{exeps} with $\tilde{\sigma}(x)=\sqrt{x}$) the terms agree with each other. Indeed for the first term this follows as (cf. \cite[Section 5]{alos2018exponentiation})
\begin{eqnarray*}
(X\diamond M)_0(T)&&=\epsilon\,\rho\,\int_0^T\int_{t_1}^T\mathbb{E}[\sqrt{u_{t_1}^\epsilon(0)}g(t_2-t_1)\sqrt{u_{t_1}^\epsilon(0)}]\,dt_2\,dt_1\\
&&=\epsilon\,\rho\,\int_0^T\mathbb{E}[u_{t_1}^\epsilon(0)]\int_{t_1}^Tg(t_2-t_1)\,dt_2\,dt_1=\epsilon\, C_0^{x,u}(u).
\end{eqnarray*}
In the generalized stochastic volatility models with affine drift given in \ref{exeps} the spot variance $u_t^\epsilon(0)$ satisfies the stochastic Volterra equation
\begin{equation*}\label{sveu}
u_t^\epsilon(0)=u(t)+\int_0^tg(t-s)\tilde{\sigma}(u_s^\epsilon(0))\,d\beta_s^1.
\end{equation*}
For this class we can show that for the terms of the forest expansion to agree with the terms of the BG expansion it is not only sufficient (as shown in \cite{alos2018exponentiation}) but also necessary that the forward variance model is affine in the sense of \cite{gatheral2018affine}.
\begin{proposition}If $u^\epsilon$ satisfies \eqref{sveu} and $(X\diamond M)_0(T)=\epsilon\, C_0^{x,u}(u)$ holds true, then $u^\epsilon$ is necessarily affine, i.e. $\tilde{\sigma}(x)=\sqrt{x}$.
\begin{proof}If $u^\epsilon$ satisfies \eqref{sveu} then $(X\diamond M)_0(T)$ is given by
\begin{eqnarray*}
&&(X\diamond M)_0(T)=\epsilon\,\rho\,\int_0^T\mathbb{E}[\tilde{\sigma}(u_{t_1}^\epsilon(0))\sqrt{u_{t_1}^\epsilon(0)}]\int_{t_1}^Tg(t_2-t_1)\,dt_2\,dt_1.
\end{eqnarray*}
Thus if $(X\diamond M)_0(T)=\epsilon\, C_0^{x,u}(u)$ then necessarily $\mathbb{E}[\tilde{\sigma}(u_{t_1}^\epsilon(0))\sqrt{u_{t_1}^\epsilon(0)}]=u(t_1)$ for all $t_1\in [0,T]$ and hence in particular for $t_1=0$ which gives $\tilde{\sigma}(u(0))=\sqrt{u(0)}$.
\end{proof}
\end{proposition}
\section{Truncated weak Taylor approximation}\label{sWTA}
In this section we will introduce the truncated version of \cite[Theorem 1]{TS2011} in a similar way \cite[Theorem 2.2]{YN1992} introduced a truncated version of \cite[Theorem 2.3]{WS1987} and end up essentially at \cite[Theorem 2.1]{KT2003} which is stated without a proof. We will make use of Malliavin calculus and use the notation of \cite{nd2006}. 
\begin{definition}\label{LND}
For a measurable subset $A\subset \Omega$ we will say that a random variable $F\in \mathbb{D}^\infty$ is \emph{locally non-degenerate on $A$} if it’s Malliavin covariance matrix $\gamma(F)$ satisfies $\gamma(F)^{-1}\in L^p$ on $A$ for all $p\geq 1$. If we can choose $A=\Omega$ we will as usual just say that $F$ is non-degenerate.
\end{definition}
The following crucial lemma can be shown just as in \cite[Proposition 2.1.4]{nd2006}. 
\begin{lemma}\label{IBP}
If $F\in \mathbb{D}^\infty$ is locally non-degenerate on $A$, then for $\phi \in C_{p}^\infty$ and $G\in \mathbb{D}^\infty$ satisfying $G=0$ on $A^c$, the integration by parts formula holds
\begin{equation*}
\mathbb{E}[\phi'(F)G]=\mathbb{E}[\phi(F)H],\quad \text{where } H:=\delta\Big(G \frac{DF}{\gamma(F)}\Big)\in \mathbb{D}^\infty.
\end{equation*}
\end{lemma}
In the following we will use the notion of Malliavin weights $\pi$ which are Skohorod integrals that will be evaluated either at $\epsilon=0$ or at an arbitrary $\epsilon \in[-1,1]$. In the former case we will just write $\pi$ and in the latter case $\pi(\epsilon)$. 
\begin{proposition}\label{TWA}
Let $\epsilon\mapsto (F_\epsilon,\eta_\epsilon)$ be a smooth map from $[-1,1]$ into $\mathbb{D}^\infty(\mathbb{R}^2)$, such that for all $\epsilon\in [-1,1]$, $F_\epsilon$ is locally non-degenerate on $\{\eta_\epsilon\in (-1,1)\}$. If $F_0$ is non-degenerate, $\eta_\epsilon$ satisfies $\eta_0=0$ and \cite[Condition 4 of Theorem 2.2]{YN1992}, i.e.
\begin{equation}\label{tcond}
\mathbb{P}\Big[\vert \eta_\epsilon\vert >\frac{1}{2}\Big]=o(\epsilon^k) \,\,\, \text{ as } \epsilon\rightarrow 0, \text{ for all } k\geq 1,
\end{equation}
then $F_\epsilon$ admits a weak Taylor approximation of arbitrary order around $0$ in the sense of \cite[Definition 2]{TS2011}, that is, for every $n\geq 1$ and $f$ in $C_p^\infty$ or bounded and Borel measurable, there are \emph{Malliavin weights} $\pi_0,...,\pi_n\in \mathbb{D}^\infty$ such that 
\begin{equation}\label{WT}
\vert \mathbb{E}[f(F_\epsilon)]-\sum_{i=0}^n\frac{\epsilon^i}{i!}\mathbb{E}[f(F_0)\pi_i]\vert=o(\epsilon^n),\quad \text{as } \epsilon \rightarrow 0.
\end{equation}
\begin{proof}
As $F_0$ is non-degenerate, we can compute the $\pi_1,\pi_2,...$ as in \cite[Theorem 1]{TS2011} but possibly without having the property \eqref{WT}. Let $\psi$ be a smooth function such that $\psi(y)=0$ if $\vert y\vert \geq 1$ and $\psi(y)=1$ if $\vert y \vert<\frac{1}{2}$. Then for any $G\in \mathbb{D}^\infty$ we have $\psi(\eta_\epsilon)G\in \mathbb{D}^\infty(\mathbb{R})$ and $\psi(\eta_\epsilon)G=0$ on $\{\vert \eta_\epsilon\vert \geq 1\}$. We show the claim first for $\phi\in C_{p}^\infty$, as in this case we can apply the integration by parts formula given in Lemma \ref{IBP} to arrive at
\begin{equation*}
\mathbb{E}[\phi'(F_\epsilon)\psi(\eta_\epsilon)G]=\mathbb{E}\Big[\phi(F_\epsilon)\delta\Big(\frac{\psi(\eta_\epsilon)GDF_\epsilon}{\gamma(F_\epsilon)}\Big)\Big].
\end{equation*}
As $\psi(\eta_0)=1$ and $\psi^{(n)}(\eta_0)=0$ for all $n\geq 1$, we can argue from this in the same way as in the proof of \cite[Theorem 1]{TS2011} to see that $\phi(F_\epsilon)\psi(\eta_\epsilon)$ admits a weak Taylor approximation of any order $n\geq 1$
\begin{equation*}\label{mainHelp}
\vert \mathbb{E}[\phi(F_\epsilon)\psi(\eta_\epsilon)]-\sum_{i=0}^n\frac{\epsilon^i}{i!}\mathbb{E}[\phi(F_0)\pi_i]\vert=o(\epsilon^n)\quad \text{ as }\,\,\epsilon \rightarrow 0.
\end{equation*}
Hence we can conclude for $\phi\in C_{p}^\infty$ by applying the H\"{o}lder inequality
\begin{eqnarray*}
&&\vert \mathbb{E}[\phi(F_\epsilon)]-\sum_{i=0}^n\frac{\epsilon^i}{i!}\mathbb{E}[\phi(F_0)\pi_i]\vert\\
&&\quad\leq \vert \mathbb{E}[\phi(F_\epsilon)]-\mathbb{E}[\phi(F_\epsilon)\psi(\eta_\epsilon)]\vert +\vert \mathbb{E}[\phi(F_\epsilon)\psi(\eta_\epsilon)]-\sum_{i=0}^n\frac{\epsilon^i}{i!}\mathbb{E}[\phi(F_0)\pi_i]\vert\\
&&\quad= \vert \mathbb{E}[\phi(F_\epsilon)(1-\psi(\eta_\epsilon))]\vert +o(\epsilon^n)\leq \sqrt{\mathbb{E}[\phi(F_\epsilon)^2]\mathbb{E}[(1-\psi(\eta_\epsilon))^2]} +o(\epsilon^n)
\end{eqnarray*}
and by recalling the definition of $\psi$, noting that $\sqrt{\mathbb{E}[\phi(F_\epsilon)^2]}\rightarrow \sqrt{\mathbb{E}[\phi(F_0)^2]}$ as $\epsilon \rightarrow 0$ and using \eqref{tcond} we see that
\begin{eqnarray*}
\sqrt{\mathbb{E}[\phi(F_\epsilon)^2]\mathbb{E}[(1-\psi(\eta_\epsilon))^2]}&=& \sqrt{\mathbb{E}[\phi(F_\epsilon)^2]\mathbb{E}[(1-\psi(\eta_\epsilon))^2\textbf{1}_{\vert \eta_\epsilon\vert >\frac{1}{2}}]}\\
&\leq&\sqrt{\mathbb{E}[\phi(F_\epsilon)^2]\mathbb{P}\Big[\Big\vert \eta_\epsilon\Big\vert >\frac{1}{2}\Big]}=o(\epsilon^{k}),
\end{eqnarray*}
for all $k\geq 1$ which gives \eqref{WT} for $\phi$ in $C_{p}^\infty$. In particular for all $n\geq 1$
\begin{equation}\label{TWAwRG}
\mathbb{E}[\phi(F_\epsilon)]=\sum_{i=0}^n\frac{\epsilon^i}{i!}\mathbb{E}[\phi(F_0)\pi_i]+\frac{\epsilon^{n+1}}{(n+1)!}\mathbb{E}[\phi(F_\xi)\psi(\eta_\xi)\pi_{n+1}(\xi)],
\end{equation}
where $\xi\in (0,\epsilon)$. To show the claim for the general case it suffices to show the validity of \eqref{TWAwRG} for all bounded and Borel measurable functions, which can be shown with the Monotone Class Theorem (cf. \cite[Chapter I, Theorem 8]{PP2005}). In fact, we note that $C_{p}^\infty$ is closed under multiplication and generates (a $\sigma$-algebra that contains) the Borel $\sigma$-algebra (as for the example the identity map on $\mathbb{R}$ is included in $C_{p}^\infty$). Now let $\mathcal{H}$ denote the set of all bounded measurable functions that satisfy \eqref{TWAwRG}. Then by the linearity of $\eqref{TWAwRG}$ in $\phi$, $\mathcal{H}$ is a vector space that trivially contains the constant functions. Now let $(f_m)_{m\geq 1}\subset \mathcal{H}$ be such that $0\leq f_1\leq f_2\leq\cdots $ and $\lim_{m\rightarrow \infty}f_m=:f$ is bounded. The claim follows from \cite[Chapter I, Theorem 8]{PP2005} if $f\in \mathcal{H}$. We show that each term in \eqref{TWAwRG} converges. For the term on the left-hand side this follows from dominated convergence. For the terms on the right-hand side this again follows from dominated convergence after an application of H\"{o}lder's inequality. In fact, recalling that the Skohorod integrals $\pi(\xi)$, for $\xi\in [0,\epsilon]$, are in $\mathbb{D}^\infty$ and as each $f_m$ as well as $f$ is bounded we can apply the H\"{o}lder's inequality and accordingly again from dominated convergence
\begin{equation*}
\vert\mathbb{E}[(f_m(F_\epsilon)-f(F_\epsilon))\pi(\epsilon)]\vert\leq \sqrt{\mathbb{E}[(f_m(F_\epsilon)-f(F_\epsilon))^2]}\sqrt{\mathbb{E}[(\pi(\epsilon))^2]}\rightarrow 0
\end{equation*}
as $m\rightarrow \infty$ and hence the claim.
\end{proof}
\end{proposition}
In the following corollary we state conditions under which the local non-degeneracy condition and \eqref{tcond} are satisfied. It corresponds to a generic version (in particular it does not explicitly depend on the first variation process) of \cite[Theorem 3.2 and Lemma A.2]{KT2003} (see also \cite[Section 5]{ty2004}). Nevertheless the proof is similar. 
\begin{corollary}\label{TWAcor}Let $F_\epsilon$ be $\mathcal{F}_t$-measurable for all $\epsilon\in [-1,1]$ such that $\epsilon\mapsto F_\epsilon$ is a smooth map from $[-1,1]$ into $\mathbb{D}^\infty$ and $F_0$ is non-degenerate such that its Malliavin covariance matrix $\gamma(F_0)$ satisfies $\gamma(F_0)^{-1}>\delta$ for some strictly positive $\delta$. If also 
\begin{equation}\label{lcond}
\sum_{i=1}^{d}\int_0^t\mathbb{E}[\vert D_s^iF_\epsilon-D_s^iF_0\vert^{2p}]\,ds=o(\epsilon^{2p}),\quad \text{ as } \epsilon \rightarrow 0,
\end{equation}
holds for all $p\geq 1$, then the conditions of Proposition \ref{TWA} are satisfied.
\begin{proof}We first construct a suitable $\eta$ and show that the local non-degeneracy of $F_\epsilon$ on $\{\eta_\epsilon\in (-1,1)\}$. Recalling that the Malliavin covariance matrix $\gamma(F_\epsilon)$ is given by $\gamma(F_\epsilon)=\sum_{i=1}^d\int_0^t( D_s^iF_\epsilon )^2\,ds$ we find
\begin{eqnarray*}
&&\vert \gamma(F_\epsilon)-\gamma(F_0)\vert\leq \sum_{i=1}^d\int_0^t\vert (D^i_sF_\epsilon)^2-(D^i_sF_0)^2\vert\,ds.
\end{eqnarray*}
By using the inequality $\vert x^2-y^2\vert \leq \vert x-y\vert^2+2\vert y\vert\vert x-y\vert$ (that holds true for all $x,y\in \mathbb{R}$) as well as the H\"{o}lder and the Cauchy-Schwarz inequality we see  
\begin{eqnarray*}
&&\vert \gamma(F_\epsilon)-\gamma(F_0)\vert\leq \sum_{i=1}^d\int_0^t\vert D^i_sF_\epsilon- D^i_sF_0\vert^2+2\vert D_s^iF_0\vert\vert D_s^iF_\epsilon-D_s^iF_0\vert\,ds\\
&&\leq\sum_{i=1}^d\int_0^t\vert D_s^iF_\epsilon - D_s^iF_0\vert^2\,ds+2\sqrt{\gamma(F_0)}\sqrt{\sum_{i=1}^d\int_0^t\vert D_s^iF_\epsilon - D_s^iF_0\vert^2\,ds}.
\end{eqnarray*}
If we now choose the map $[-1,1]\ni\epsilon\mapsto \eta_\epsilon\in \mathbb{D}^\infty$ as 
\begin{equation*}
\eta_\epsilon=c\,\frac{\sum_{i=1}^{d}\int_0^t\vert D_s^iF_\epsilon-D_s^iF_0\vert^2\,ds}{\gamma(F_0)},
\end{equation*}
for some $c>0$ satisfying $\delta:=1-(\frac{1}{c}+\frac{2}{\sqrt{c}})>0$, we see that on $\{\eta_\epsilon<1\}$ 
\begin{eqnarray*}\label{etaeps2}
&&\vert \gamma(F_\epsilon)-\gamma(F_0)\vert\leq \Big(\frac{1}{c}+\frac{2}{\sqrt{c}}\Big)\gamma(F_0)
\end{eqnarray*}
and accordingly 
\begin{eqnarray*}
\gamma(F_\epsilon)&=&\gamma(F_0)+\gamma(F_\epsilon)-\gamma(F_0)\geq \gamma(F_0)-\vert\gamma(F_\epsilon)-\gamma(F_0)\vert\nonumber\\
&\geq&\gamma(F_0)\Big(1-\Big(\frac{1}{c}+\frac{2}{\sqrt{c}}\Big)\Big)= \delta\gamma(F_0)\label{etaeps3},
\end{eqnarray*}
which gives the local non-degeneracy of $F_\epsilon$ on the set $\{\eta_\epsilon<1\}$. Hence the claim follows if Property \eqref{tcond} is satisfied, which it is, as 
\begin{eqnarray*}
P[\eta_{\epsilon}>\frac{1}{2}]&&\leq 2^p \mathbb{E}[\eta_\epsilon^p]= 2^p \mathbb{E}[\gamma(F_0)^{-p}c^p(\sum_{i=1}^{d}\int_0^t\vert D_s^iF_\epsilon-D_s^iF_0\vert^2\,ds)^p]\\
&&\leq 2^p \delta^p C(t) \sum_{i=1}^{d}\int_0^t\mathbb{E}[\vert D_s^iF_\epsilon-D_s^iF_0\vert^{2p}\,ds]=o(\epsilon^{2p}),
\end{eqnarray*}
where we applied the Chebyshev inequality, the Jensen inequality and \eqref{lcond}.
\end{proof}
\end{corollary}
\section{Vol-of-vol Expansion: the regular case}\label{SNE}
The \emph{vol-of-vol expansion} up to order $n$ holds true, if for every $f$ that is either in $C_p^\infty$ or bounded and measurable there are \emph{Malliavin weights} $\pi_i\in \mathbb{D}^\infty$ such that 
\begin{equation}\label{WTX}
\vert \mathbb{E}[f(X_T^\epsilon)]-\sum_{i=0}^n\frac{\epsilon^i}{i!}\mathbb{E}[f(X_T^0)\pi_i]\vert=o(\epsilon^n),\quad \text{as } \epsilon \rightarrow 0.
\end{equation}
In this section we look at the regular version of this expansion, i.e. to the case where in \eqref{ForwV} the function $g$ belongs to $H$. We will write for notational convenience $\Sigma(u):=m(\sigma(u),g)$ and notice that under the given conditions $\Sigma\in C_b^\infty(H)$. We assume further that the square root function appearing in the log-price \eqref{logPrice} is replaced by a smooth approximation $w\in C^\infty_{bb}$. In this case the following holds.
\begin{lemma}\label{maliTHM}Under the conditions given above, $(X_t^\epsilon,u_t^\epsilon)\in \mathbb{D}^\infty(\mathbb{R}\times H)$ for each $\epsilon\in [-1,1]$. In particular, for $r\in [0,t]$, we have the following representations
\begin{eqnarray*}
D_r^1X_t^\epsilon&&=w(u_r^\epsilon(0))\rho-\frac{1}{2}\int_r^tD_r^1u_{t_1}^\epsilon(0)\,dt_1\\
&&+\rho\int_r^t dw(u_{t_1}^\epsilon(0))D_r^1u_{t_1}^\epsilon(0)\,d\beta^1_{t_1}\label{DX1},\\
D_r^2X_t^\epsilon&&=w(u_r^\epsilon(0))\sqrt{1-\rho^2}\label{DX2},\\
D_r^1u_t^\epsilon&&=\epsilon S_{t-r}\Sigma(u_{r}^\epsilon)+\epsilon\int_r^tS_{t-t_1}d\Sigma(u_{t_1}^\epsilon)D_r^1u_{t_1}^\epsilon\,d\beta_{t_1}^1\label{Du}.
\end{eqnarray*}
\begin{proof}
See for example \cite[Lemma 5.3]{ct2007} or \cite[Theorem IV 2.17]{OA2016}.
\end{proof}
\end{lemma}
\begin{proposition}\label{propSNE}Under the conditions given above, for every strictly positive initial curve $u$, $X^\epsilon$ admits a truncated weak Taylor expansion of arbitrary order.
\begin{proof}We will show that the conditions of Corollary \ref{TWAcor} are satisfied. From the strict positivity of $u$ it follows that $X_t^0$ is non-degenerate, as its Malliavin matrix 
\begin{equation*}
\gamma(X_t^0)=\sum_{i=1}^2\int_0^t(D^i_sX_t^0)^2\,ds=\int_0^tw(u(s))^2\,ds
\end{equation*}
is deterministic and due the given condition strictly positive. It remains to show that \eqref{lcond} holds true, i.e.
\begin{equation}\label{SNEhelp3}
\sum_{i=1}^2\int_0^t\mathbb{E}[\vert D_s^iX_t^\epsilon-D_s^iX_t^0\vert^{2p}]\,ds=o(\epsilon^{2p})\quad \text{as }\epsilon \rightarrow 0.
\end{equation}
Recalling Lemma \ref{maliTHM} and using the Jensen inequality gives
\begin{eqnarray*}
\mathbb{E}[\vert D_s^1X_t^\epsilon-D_s^1X_t^0\vert^{2p}]&&\leq 3^{2p-1}\rho^{2p}\mathbb{E}[\vert w(u_s^\epsilon(0))-w(u_s^0(0))\vert^{2p}]\\
&&+3^{2p-1}\frac{1}{2^{2p}}\mathbb{E}[\vert \int_s^tD_s^1u_r^\epsilon(0)\,dr \vert^{2p}]\\
&&+3^{2p-1}\mathbb{E}[\vert \int_s^tdw(u_r^\epsilon(0))D_s^1u_r^{\epsilon}(0)\,d\beta_r^1\vert^{2p}].
\end{eqnarray*}
By recalling $w\in C_{bb}^\infty$ (hence in particular Lipschitz continuous) and that the point evaluation in $H$ is continuous we see that by applying again the Jensen inequality to the second term and the Burkholder-Davis-Gundy inequality to the third term, we arrive at the inequality
\begin{equation}\label{SNEhelper5}
\sum_{i=1}^2\mathbb{E}[\vert D_s^iX_t^\epsilon-D_s^iX_t^0\vert^{2p}]\leq C_1 \mathbb{E}[\Vert u_s^\epsilon-u_s^0\Vert^{2p}]+C_2 \int_s^t\mathbb{E}[\Vert D_s^1u_r^\epsilon\Vert^{2p} ]\,dr
\end{equation}
for positive constants $C_1,C_2$ that depend only on $t$. In the following we will introduce further constants $C_3-C_{10}$ that share this properties. We need to show that both terms in \eqref{SNEhelper5} are in $o(\epsilon^{2p})$ when $\epsilon \rightarrow 0$ as in this case \eqref{SNEhelp3} will follow. For the first term this follows from \cite[Theorem 7.2]{da2014stochastic} by noting that $u_s^\epsilon-u_s^0=\epsilon u_s^1$ with $u_s^1$ starting at $0$, hence
\begin{equation*}
\sup_{s\in[0,t]}\mathbb{E}[\Vert u_s^\epsilon-u_s^0\Vert^{2p}]=\epsilon^{2p}\sup_{s\in[0,t]}\mathbb{E}[\Vert u_s^1\Vert^{2p}]\leq \epsilon^{2p} C_3=o(\epsilon^{2p}).
\end{equation*}
The claim follows upon showing that the second term in \eqref{SNEhelper5} satisfies
\begin{equation}\label{SNEhelper7}
\sup_{s\in [0,t]}\int_s^t\mathbb{E}[\Vert D_s^1u_r^\epsilon\Vert^{2p}]\,dr=o(\epsilon^{2p}).
\end{equation}
By recalling Lemma \ref{maliTHM}, $\Sigma\in C_b^\infty$ (and hence the boundedness of $d\Sigma$ and linear growth property of $\Sigma$), that there is a $M\geq 1$ and a real number $\omega$ such that $\Vert S_t\Vert\leq Me^{\omega t}$ as $\{S_t\,\vert\,t\geq 0\}$ is a strongly continuous semigroup, then applying the Jensen inequality and \cite[Proposition 7.3 and Theorem 7.2]{da2014stochastic} we arrive at
\begin{eqnarray*}
&&\sup_{r\in [s,t]}\mathbb{E}[\Vert D_s^1u_r^\epsilon\Vert^{2p}]\\
&&\quad\quad\leq\epsilon^{2p}C_{4}\sup_{r\in [s,t]}\mathbb{E}[\Vert S_{r-s}\Sigma(u_s^\epsilon)\Vert^{2p}]+\epsilon^{2p}C_{5}\int_s^t\mathbb{E}[\Vert D_s^1 u_l^\epsilon\Vert^{2p}]\,dl\\
&&\quad\quad\leq \epsilon^{2p}C_{6}\sup_{r\in [s,t]}\mathbb{E}[\Vert u_s^\epsilon\Vert^{2p}]+\epsilon^{2p}C_{7}\int_s^t\sup_{r\in [s,l]}\mathbb{E}[\Vert D_s^1u_r^\epsilon\Vert^{2p}]\,dl\\
&&\quad\quad\leq\epsilon^{2p}C_{8}(1+\Vert u\Vert^{2p})+\epsilon^{2p}C_{7}\int_s^t\sup_{r\in [s,l]}\mathbb{E}[\Vert D_s^1u_r^\epsilon\Vert^{2p}]\,dl.
\end{eqnarray*}
Finally, applying the Gronwall inequality
\begin{eqnarray*}
&&\sup_{r\in [s,t]}\mathbb{E}[\Vert D_s^1u_r^\epsilon\Vert^{2p}]\leq \epsilon^{2p}C_9e^{\epsilon^{2p}C_{10}}
\end{eqnarray*}
gives Property \eqref{SNEhelper7} and thus the claim.
\end{proof}
\end{proposition}
\begin{remark}[Approximation of the square root]\label{sqremark}
It is shown in \cite[Section 4.1]{TT2012} that the impact of the approximation of the square root in the context of small noise expansions is asymptotically negligible and hence we will mostly omit mentioning it.
\end{remark}
\section{Push-down Malliavin weights}\label{PDMW}
The Malliavin weights arise by iteratively applying Malliavin's integration by parts formula
\begin{equation*}\label{push1}
\frac{\partial^n}{\partial \epsilon^n}\Big\vert_{\epsilon=0}\mathbb{E}[f(X_T^\epsilon)]=\mathbb{E}[f(X_T^0)\pi_n]
\end{equation*}
where $\pi_n:=\pi(0)$. Then $\pi_0(\epsilon):=1$ and for $n\geq 1$
\begin{equation*}\label{push2}
\pi_n(\epsilon):=\delta(X_T^{(1),\epsilon}\frac{DX_T^\epsilon}{\Vert DX_T^\epsilon\Vert^2}\pi_{n-1}(\epsilon))+\frac{\partial}{\partial \epsilon}\pi_{n-1}(\epsilon),
\end{equation*}
where $X_T^{(1),\epsilon}:=\frac{\partial}{\partial \epsilon}X_T^{\epsilon}$. A version of the Malliavin weights that is particularly useful are the so-called \emph{push-down Malliavin weights} (we adopted this name from \cite{TT2012}) which are given by $\mathbb{E}[\pi_n\vert Y_T]$ (recall \eqref{YSigma}) for $n\geq 0$ and satisfy
\begin{equation}\label{push3}
\mathbb{E}[f(X_T^0)\pi_n]=\mathbb{E}[f(X_T^0)\mathbb{E}[\pi_n\vert Y_T]],\quad \text{ for } n\geq 0.
\end{equation}
We will now derive a representation of the expectations in \eqref{push3} that is free from the recursively defined Malliavin weights (see above and compare also to Representations (53) and (54) in \cite{TT2012}). We need the following results related to the Hermite and Bell polynomials. 
\begin{lemma}\label{Hrep}Let $X\in L^2(\Omega)$ and $Z$ be a real normally distributed random variable with zero mean and variance given by $\Sigma$. Then
\begin{equation*}\label{ans}
\mathbb{E}[X\,\vert\,Z]=\sum_{n=0}^\infty c_n H_n(Z,\Sigma), \text{ where }\quad c_n:=\frac{n!}{\Sigma^{n}}\mathbb{E}[XH_n(Z,\Sigma)].
\end{equation*}
\begin{proof}
Cf. \cite[Lemma 1]{ttt2012}.
\end{proof}
\end{lemma}
\begin{lemma}\label{MBS}Let MBs be the function that gives the maximum number of Brownian motions appearing in one term of the sum given in \eqref{Bnp}. Then, 
\begin{equation}\label{mbkn}
MBs\Big(B_{n,k}(X_T^{(1)},X_T^{(2)},...,X_T^{(n-k+1)})\Big)=k+n
\end{equation}
and accordingly 
\begin{equation}\label{corpeq}
\mathbb{E}[B_{n,k}(X_T^{(1)},...,X_T^{(n-k+1)})H_l(T)]=0\quad \text{for } l>n+k.
\end{equation}
\begin{proof}
The equality \eqref{mbkn} follows from It\^{o}'s product rule, i.e. 
\begin{eqnarray}\label{um}
&&u_T^{(m)}(0)=m\int_0^T\frac{\partial^{m-1}}{\partial \epsilon^{m-1}}\Big\vert_{\epsilon=0}\Sigma(u_{t_1}^\epsilon)(T-t_1)\,d\beta_{t_1}^1
\end{eqnarray}
gives $MBs(u_T^{(m)}(0))=m$ and $MBs(u_T^{(m)}(0)u_T^{(n)}(0))=m+n$ from the product rule. Accordingly from the chain rule $MBs(\sqrt{u_T(0)}^{(m)})=m$. Then from
\begin{equation}\label{xmjm}
\begin{cases}
&(X_T^{(m)})^{j_m}=\int_0^Tj_m(X_{t_1}^{(m)})^{j_m-1}\sqrt{u_{t_1}(0)}^{(m)}\,d\beta_{t_1}\\
&+\frac{1}{2}\int_0^Tj_m(j_m-1)(X_{t_1}^{(m)})^{j_m-2}(\sqrt{u_{t_1}(0)}^{(m)})^2-j_m(X_{t_1}^{(m)})^{j_m-1}u_{t_1}^{(m)}(0)\,dt_1
\end{cases}
\end{equation}
with $j_m=1$ we see that $MBs(X_T^{(m)})=1+MBs(\sqrt{u_T(0)}^{(m)})=1+m$ and \eqref{mbkn} and hence \eqref{corpeq} follows from another application of the product rule and \eqref{tnk} as
\begin{equation*}
MBs\Big(\prod_{m=1}^{n-k+1} (X_T^{(m)})^{j_m}\Big)=\sum_{m=1}^{n-k+1} MBs( (X_T^{(m)})^{j_m})=\sum_{m=1}^{n-k+1} j_m(1+m).
\end{equation*}
\end{proof}
\end{lemma}
\begin{lemma}\label{Hrep2}The following equality holds true
\begin{equation*}
\frac{d^n}{dx^n}\mathbb{E}[f(X_T^0)]=\frac{n!}{\Sigma_T^{n}}\mathbb{E}[f(X_T^0)H_n(Y_T, \Sigma_T)], \quad n\geq 0.
\end{equation*}
\begin{proof}This is straight forward by recalling $X_T^0=x+Y_T-\frac{1}{2}\Sigma_T=Z_T-\frac{1}{2}\Sigma_T$ with $Z_T$ being normally distributed with mean $x$ and variance $\Sigma_T$. Then
\begin{eqnarray*}
&&\mathbb{E}[f(X_T^0)\frac{n!}{\Sigma_T^n}H_n(Y_T,\Sigma_T)]\\
&&=\frac{1}{\sqrt{2\pi\Sigma_T}}\int_{\mathbb{R}}f(z-\frac{1}{2}\Sigma_T)e^{-\frac{(z-x)^2}{2\Sigma_T}}(-(1)^n e^{\frac{(z-x)^2}{2\Sigma_T}}\frac{d^n}{dy^n}\Big\vert_{y=z-x}e^{\frac{-y^2}{2\Sigma_T}})\,dz\\
&&=\frac{1}{\sqrt{2\pi\Sigma_T}}\int_{\mathbb{R}}f(z-\frac{1}{2}\Sigma_T)(-(1)^n)\frac{d^n}{dy^n}\Big\vert_{y=z-x}e^{\frac{-y^2}{2\Sigma_T}}\,dz\\
&&=\frac{1}{\sqrt{2\pi\Sigma_T}}\int_{\mathbb{R}}f(z-\frac{1}{2}\Sigma_T)\frac{d^n}{dx^n}e^{\frac{-(x-z)^2}{2\Sigma_T}}\,dz=\frac{d^n}{dx^n}\mathbb{E}[f(X_T^0)].
\end{eqnarray*}
\end{proof}
\end{lemma}
Now we can return to the mentioned explicit representation of the expectations given in \eqref{push3}. We recall the Fa\`{a} di Brunno's formula in terms of Bell polynomials
\begin{equation}\label{pidual}
\frac{\partial^n}{\partial \epsilon^n}\Big\vert_{\epsilon=0}\mathbb{E}[f(X_T^{\epsilon})]=\sum_{k=1}^n\mathbb{E}[f^{(k)}(X_T^0)B_{n,k}(X_T^{(1)},...,X_T^{(n-k+1)})].
\end{equation}
\begin{proposition}\label{expREep}For $n\geq 0$ and $f$ either in $C_p^\infty$ or bounded, we have
\begin{equation*}\label{ma1}
\mathbb{E}[f(X_T^0)\pi_n]=\sum_{k=1}^n\sum_{l=k}^{n+2k} \Big(\frac{d^l}{dx^l}\mathbb{E}[f(X_T^0)]\Big)\mathbb{E}\Big[B_{n,k}(X_T^{(1)},...,X_T^{(n-k+1)}) H_{l-k}(T)\Big].
\end{equation*}
\begin{proof}
It follows from the Lemmas \ref{Hrep} and \ref{Hrep2}, that
\begin{equation}\label{ma2}
\mathbb{E}[f(X_T^0)\pi_n]=\sum_{l=0}^\infty \frac{d^l}{dx^l}\mathbb{E}[f(X_T^0)]\mathbb{E}[\pi_n H_l(T)].
\end{equation}
Now from \eqref{pidual} we see that for every $f$ that is in $C^\infty$ 
\begin{equation*}
\mathbb{E}[f(X_T^0)\pi_n]=\sum_{k=1}^n\mathbb{E}[B_{n,k}(X_T^{(1)},...,X_T^{(n-k+1)})d^k f(X_T^0)],
\end{equation*}
holds true. As $H_l(\cdot-x-\frac{1}{2}\Sigma_T,\Sigma_T)\in C^\infty_p$ and recalling \eqref{Hn1} it follows that
\begin{equation*}
\mathbb{E}[H_l(T)\pi_n]=\sum_{k=1}^n\mathbb{E}[B_{n,k}(X_T^{(1)},...,X_T^{(n-k+1)})H_{l-k}(T)].
\end{equation*}
Plugging this into \eqref{ma2} gives
\begin{equation*}
\mathbb{E}[f(X_T^0)\pi_n]=\sum_{k=1}^n\sum_{l=k}^\infty \frac{d^l}{dx^l}\mathbb{E}[f(X_T^0)]\mathbb{E}\Big[B_{n,k}(X_T^{(1)},...,X_T^{(n-k+1)}) H_{l-k}(T)\Big]
\end{equation*}
and hence after recalling \eqref{corpeq} from Lemma \ref{MBS} the claim.
\end{proof}
\end{proposition}
\section{Vol-of-vol Expansion: the general case}\label{RFVM}
We show now the validity of the expansion in the general case where $g\in L^2([0,T])$. We will do so by approximating the right-hand side of \eqref{WTX}, uniformly in $\epsilon\in [-1,1]$, and then conclude with a change of limits argument. 
\begin{lemma}\label{hdens}Let $T\in (0,\infty)$. Then $H$ is dense in $L^2([0,T])$.
\begin{proof}
Let $\mathcal{P}$ denote the set of polynomials on $\mathbb{R}_+$. As $\mathcal{P}$ restricted to $[0,T]$ are dense in $L^2([0,T])$, it is sufficient to show that for every $P\in \mathcal{P}$ we can find a function $\tilde{P}\in H$ such that $\tilde{P}=P$ on $[0,T]$. In fact, $\tilde{P}:=m(P,f)$ satisfies this conditions, where $f$ is a smooth function such that for $0<T<T'$, it satisfies $f(x)=1$ for $x\in [0,T]$ and $f(x)=0$ for $x\geq T'$. 
\end{proof}
\end{lemma}
We denote the forward variance in \eqref{ForwVmild} for $g=g_n$ by $u^{\epsilon,n}$ and the corresponding log-price by $X^{\epsilon,n}$. The following is known for the case where $g$ corresponds to the power kernel of the rough models and $(g_n)$ to the so-called Markovian lifts (cf. \cite[Theorem 1. (b)]{harms2019strong} and \cite[Theorem A.2]{abi2019lifting}).
\begin{lemma}\label{Xg}Let $f$ be Lipschitz continuous and $(g_n)\subset H$ such that $g_n\rightarrow g$ in $L^2([0,T])$, then $\vert\mathbb{E}[f(X_T^\epsilon)-f(X_T^{\epsilon,n})]\vert\rightarrow 0$ uniformly in $\epsilon\in [-1,1]$.
\begin{proof}
It suffices to show for an arbitrary $\epsilon \in [-1,1]$ (for $\epsilon =0$ this is trivial) $\mathbb{E}[X_T^{\epsilon,n}]$ is a Cauchy sequence. We will again approximate the square root with a $C^\infty_{bb}$ function $w$. Then from the Burkholder-Gundy inequality, the Jensen inequality, the Lipschitz continuity of $w$ as well as the (local) martingality of $u_{t_1}^{\epsilon,n}(0)$ we have
\begin{eqnarray*}
&&\mathbb{E}[\vert X_T^{\epsilon,n}-X_T^{\epsilon,m}\vert^2]\leq \frac{1}{2}\mathbb{E}[\vert \int_0^T u_{t_1}^{\epsilon,n}(0)-u_{t_1}^{\epsilon,m}(0)\,dt_1\vert^2]\\
&&\quad\quad+2\mathbb{E}[\vert\int_0^T w(u_{t_1}^{\epsilon,n}(0))-w(u_{t_1}^{\epsilon,m}(0))\,d\beta_{t_1}\vert^2]\\
&&\quad\quad\leq C(T) \int_0^T \mathbb{E}[\vert u_{t_1}^{\epsilon,n}(0)-u_{t_1}^{\epsilon,m}(0)\vert^2]\,dt_1,
\end{eqnarray*}
and see that the claim follows if $\mathbb{E}[ \vert u_{t_1}^{\epsilon,n}(0)-u_{t_1}^{\epsilon,m}(0)\vert^2]$ vanishes as $n,m\rightarrow \infty$. From the Jensen inequality and the Lipschitz continuity of $\sigma$ we have
\begin{eqnarray*}
&&\mathbb{E}[ \vert u_{t_1}^{\epsilon,n}(0)-u_{t_1}^{\epsilon,m}(0)\vert^2]\\
&&\quad\quad\leq C_1 \mathbb{E}[ \int_0^{t_1} \Vert S_{t_1-t_2}\Vert^2\vert (\sigma(u_{t_2}^{\epsilon,n}(0))-\sigma(u_{t_2}^{\epsilon,m}(0))\vert^2\vert g_n(t_1-t_2)\vert^2\,dt_2]\\
&&\quad\quad+ C_1 t_1 \mathbb{E}[ \int_0^{t_1}\Vert S_{t_1-t_2}\Vert^2 \vert\sigma(u_{t_2}^{\epsilon,m}(0))\vert^2\vert g_n(t_1-t_2)-g_m(t_1-t_2)\vert^2\,dt_2]\\
&&\quad\quad\leq C_2(t_1) \int_0^{t_1}\mathbb{E}[\vert u_{t_2}^{\epsilon,n}(0)-u_{t_2}^{\epsilon,m}(0)\vert^2]\vert g_n(t_1-t_2)\vert^2\,dt_2\\
&&\quad\quad+C_2(t_1) \int_0^{t_1}\vert(g_n(t_1-t_2)-g_m(t_1-t_2)\vert^2\,dt_2,
\end{eqnarray*}
where we have used the boundedness of $\sigma$, and that $\Vert S_{t_1-t_2} \Vert\leq M e^{\omega (t_1-t_2)}$ for a $M\geq 1$ and a real $\omega$. Thus if $\omega >0$ then $\Vert S_{t_1-t_2} \Vert\leq M e^{\omega t_1}$ and else $\Vert S_{t_1-t_2} \Vert\leq M $, i.e. $C_2(t_1)$ is either constant or strictly increasing in $t_1$ and hence, as the term $ \int_0^{t_1}\vert(g_n(t_1-t_2)-g_m(t_1-t_2)\vert^2\,dt_2$ is strictly increasing with $t_1$ the claim follows from the following version of the Gronwall inequality 
\begin{eqnarray*}
&&\mathbb{E}[ \vert u_{t_1}^{\epsilon,n}(0)-u_{t_1}^{\epsilon,m}(0)\vert^2]\leq C_5(t_1) \exp\Big(\int_0^{t_1}\vert g_n(t_1-t_2)\vert^2\,dt_2\Big)\\
&&\quad\quad\quad\quad\quad\quad\quad\quad\quad\quad\quad\quad\quad\quad \int_0^{t_1}\vert g_n(t_1-t_2)-g_m(t_1-t_2)\vert^2\,dt_2.
\end{eqnarray*}
\end{proof}
\end{lemma}

For real-valued continuous semi-martingales $Y^1,...,Y^{n-k+1}$ we introduce a convenient representation of $\prod_{i=1}^{n-k+1}(\frac{Y_T^{i}}{i!})^{j_i}$, for a given $\mathbf{j}:= (j_1,...,j_{n-k+1})$ satisfying $\sum_{i=1}^{n-k+1}j_i=k$ and $\sum_{i=1}^{n-k+1}ij_i=n$. For $s=(s^1,s^2)\in \mathbb{N}\times \mathbb{N}_0$ we set $Y_T^{s}:=\langle Y^{s^1},Y^{s^2}\rangle_T$ if $s^2\neq 0$ and else $Y_T^{s}:=Y_T^{s^1}$. We define
\begin{eqnarray*}
Q_{k}^i:=&&\Big\{s=(s_1,...,s_{k-i})\in(\mathbb{N}\times \mathbb{N}_0)^{k-i}\,\Big\vert\, \sum_{l=1}^{k-i}1_{s_l^2\neq 0}=i\Big\},
\end{eqnarray*}
for $i=0,...,\lfloor \frac{k}{2}\rfloor$. By setting $\sum_{(s^1,s^2)=1}^{n-k+1}:=\sum_{s^1=1}^{n-k+1}\sum_{s^2=1}^{n-k+1}$ if $s^2\neq 0$ and else $\sum_{(s^1,s^2)=1}^{n-k+1}:=\sum_{s^1=1}^{n-k+1}$ as well as $j_0:=1$ and $j_s:=j_{s^1}j_{s^2}$ we can define
\[\sum_{n,k,s,\mathbf{j}}^i:=\sum_{s\in Q_{k}^i}\sum_{s_{1}=1}^{n-k+1}\cdots \sum_{s_{k-i}=1}^{n-k+1}\frac{\prod_{l=1}^{k-i}j_{s_l}}{2^i}\]
and arrive after $k$ applications of It\^{o}'s formula at the representation
\begin{eqnarray}\label{prodrep}
\prod_{i=1}^{n-k+1}(Y_T^i)^{j_i}&&=\sum_{i=0}^{\lfloor \frac{k}{2}\rfloor}\sum_{n,k,s,\mathbf{j}}^i \int_0^T\int_0^{t_1}\cdots \int_0^{t_{k-i-1}}\,dY_{t_{k-i}}^{s_{k-i}}\cdots \,dY_{t_1}^{s_1}.
\end{eqnarray}
\begin{lemma}\label{GGN}Let $(g_n)$ be a sequence in $H$ such that $g_n\rightarrow g$ as $n\rightarrow \infty$ in $L^2([0,T])$. Then for fixed but arbitrary $p\geq 1$ and $(k,l)\in \{1,...,p\}\times \{k,...,p+2k\}$, $\vert \mathbb{G}_{k,l}(g_n)(T;p)-\mathbb{G}_{k,l}(\kappa)(T;p) \vert$ vanishes as $n\rightarrow \infty$, where
\begin{equation}\label{Ggn}
\mathbb{G}_{k,l}(g_n)(T;p):=\mathbb{E}[B_{p,k}(X_T^{(1),n},...,X_T^{(p-k+1),n})H_{l-k}(T)].
\end{equation}
\begin{proof}
Recalling \eqref{partBell} and \eqref{tnk} we see that
\begin{eqnarray*}
&&B_{p,k}(X_T^{(1)},...,X_T^{(p-k+1)})=\sum_{\mathbf{j}\in T(p,k)}\frac{p!\prod_{i=1}^{p-k+1}\frac{1}{(i!)^{j_i}}}{\prod_{i=1}^{p-k+1}j_i!}\prod_{i=1}^{p-k+1}(X_T^{(i)})^{j_i}
\end{eqnarray*}
and applying \eqref{prodrep} to $\prod_{i=1}^{p-k+1}(X_T^{(i)})^{j_i}$ gives
\begin{eqnarray}
&&\prod_{i=1}^{p-k+1}(X_T^{(i)})^{j_i}=\sum_{i=0}^{\lfloor \frac{k}{2}\rfloor}\sum_{p,k,s,\mathbf{j}}^i \int_0^T\int_0^{t_1}\cdots \int_0^{t_{k-i-1}}\,dX_{t_{k-i}}^{(s_{k-i})}\cdots \,dX_{t_1}^{(s_1)}\label{Bnp},
\end{eqnarray}
where for $q\in \{1,...,k-i\}$
\begin{equation*}
X_{t}^{(s_q)}=
\begin{cases}
&\int_0^{t}\sqrt{u_{t_{1}}(0)}^{(s_q^1)}\,d\beta_{t_{1}}^1-\frac{1}{2}\int_0^{t}u_{t_{1}}^{(s_q^1)}(0)\,dt_{1},\quad \text{if } s_q^2=0,\\
&\int_0^{t}\sqrt{u_{t_{1}}(0)}^{(s_q^1)}\sqrt{u_{t_{1}}(0)}^{(s_q^2)}\,dt_{1},\quad \text{if } s_q^2\neq 0.
\end{cases}
\end{equation*}
The claim follows, if for arbitrary but fixed $k\in \{1,...,p\}$, $i\in \{0,...,\lfloor \frac{k}{2}\rfloor\}$, $s\in Q_{k}^i$ and for all $l\in \{0,...,p+k\}$
\begin{equation*}
\begin{cases}
&\mathbb{E}[\int_0^T\int_0^{t_1}\cdots \int_0^{t_{k-i-1}}\,dX_{t_{k-i}}^{(s_{k-i})}\cdots \,dX_{t_1}^{(s_1)}H_l(T)]\\
&\quad \rightarrow \mathbb{E}[\int_0^T\int_0^{t_1}\cdots \int_0^{t_{k-i-1}}\,dX_{t_{k-i}}^{(s_{k-i}),n}\cdots \,dX_{t_1}^{(s_1,n)}H_l(T)]\quad \text{as } n\rightarrow \infty,
\end{cases}
\end{equation*}
holds true. In the following we set $k_i:=k-i$ for notational convenience.

As the initial curve $u$ is bounded and $\sigma\in C_b^\infty(H)$ it follows from the product structure $m(\sigma(u),g)$ that it means no loss of generality to assume that $u\equiv 1$ and $\sigma(u)\equiv u$. Further, as in this case
\begin{eqnarray*}
\sqrt{u_t(x)}^{(m)}&&=\frac{1}{2}u^{(m)}_t(x)+\sum_{l=2}^m \sqrt{1}^{(l)}B_{m,l}(u^{(1)}_t(x),...,u^{(m-l+1)}_t(x)),
\end{eqnarray*}
it means no loss of generality when it comes to proving the claim, to look at $\tilde{X}$ instead of $X$, where
\begin{equation*}
\tilde{X}_{t}^{(s_q)}:=
\begin{cases}
&\int_0^{t_1}u^{(s_q)}_{t_1}(0)\,d\beta_{t_{1}}^1,\quad \text{if } s_q^2=0,\\
&\int_0^{t_{1}}u^{(s_q)}_{t_1}(0)\,dt_{1},\quad \text{if } s_q^2\neq 0,
\end{cases}
\end{equation*}
where $u^{(s_q)}_{t}(x):=u^{(s_q^1)}_{t}(x)$ if $s_q^2=0$ and else $u^{(s_q)}_{t}(x):=u^{(s_q^1)}_{t}(x)u^{(s_q^2)}_{t}(x)$. We define $f_{j}(t):=1$ for $j=k_i$ and for $j=0,...,k_i-1$
\begin{eqnarray*}
f_{j}(t)&&:=\int_0^t f_{j+1}(t_1)\,d\tilde{X}_{t_1}^{(s_{j+1})}=\int_0^t\int_0^{t_1}\cdots\int_0^{t_{k_i-1-j}}\,d\tilde{X}_{t_{k_i-j}}^{(s_{k_i})}\cdots\,d\tilde{X}_{t_1}^{(s_{j+1})}.
\end{eqnarray*}
With this notation, the claim follows, if for  $l=0,...,k+p$ 
\begin{eqnarray}\label{induclaim}
\mathbb{E}[(f_{0}(T)-f^n_{0}(T))H_l(T)] \rightarrow 0 \quad \text{as } n\rightarrow 0.
\end{eqnarray} 
We will assume that $s^2_1,...,s^2_i\neq 0$ and hence $s_{i+1}^2,...,s_{k-i}^2=0$. Then
\begin{equation*}
f_{j}(t_j)=
\begin{cases}
&\int_0^{t_j} f_{j+1}(t_{j+1})u^{(s_{j+1})}_{t_{j+1}}(0)\,dt_{j+1},\quad \text{for }j=0,...,i-1,\\
&\int_0^{t_j} f_{j+1}(t_{j+1})u^{(s_{j+1})}_{t_{j+1}}(0)\,d\beta_{t_{j+1}},\quad\text{for }j=i,...,k_i-1.
\end{cases}
\end{equation*}
We now derive a representation of a generic term of $\mathbb{E}[f_{0}(T)H_l(T)]$, i.e.
\begin{eqnarray*}
&&\mathbb{E}[f_{0}(T)H_l(T)] =\int_0^T\mathbb{E}[f_1(t_1)u_{t_1}^{(s_1)}(0)H_l(t_1)]\,dt_1\\
&&\quad\quad=\int_0^T\int_0^{t_1}\mathbb{E}[f_2(t_2)u_{t_2}^{(s_2)}(0)u_{t_2}^{(s_1)}(0)H_l(t_2)]\,dt_2\,dt_1+...\\
&&\quad\quad=\int_0^T\int_0^{t_1}\cdots\int_0^{t_{i-1}}\mathbb{E}[f_i(t_i)\prod_{j=1}^iu_{t_i}^{(s_j)}(0)H_l(t_i)]\,dt_i\cdots\,dt_1+... .
\end{eqnarray*} 
As $f_i(t_i)=\int_0^{t_i}\cdots\int_0^{t_{k_i-1}}u_{t_{k_i}}^{(s_{k_i})}(0)\,d\beta_{t_{k_i}}^1\cdots u_{t_{i+1}}^{(s_{i+1})}(0)\,d\beta_{t_{i+1}}^1$, we can continue with (again choosing one generic term)
\begin{eqnarray*}
&&\mathbb{E}[f_i(t_i)\prod_{j=1}^iu_{t_i}^{(s_j)}(0)H_l(t_i)]\\
&&\quad \quad\quad=\int_0^{t_i}\mathbb{E}[f_{i+1}(t_{i+1})\prod_{j=1}^{i+1}u_{t_{i+1}}^{(s_j)}(0)H_{l-1}(t_{i+1})]\,d_{t_{i+1}}+...\\
&&\quad \quad\quad=\int_0^{t_i}\int_0^{t_{i+1}}\mathbb{E}[f_{i+2}(t_{i+2})\prod_{j=1}^{i+2}u_{t_{i+2}}^{(s_j)}(0)H_{l-2}(t_{i+2})]d_{t_{i+2}}d_{t_{i+1}}+...\\
&&\quad \quad\quad=\int_0^{t_i}\cdots\int_0^{t_{k_i-1}}\mathbb{E}[\prod_{j=1}^{k_i}u_{t_{k_i}}^{(s_j)}(0)H_{l-(k-2)}(t_{k_i})]d_{t_{k_i}}\cdots\,d_{t_{i+1}}+...
\end{eqnarray*}
We recall that $u_t^{(0)}(0)=1$ and for $m\in \mathbb{N}$ 
\begin{eqnarray*}
u_{t_q}^{(m)}(0)&&=m!\int_0^{t_q}\cdots \int_0^{t_{q+m-1}}g(t_q-t_{q+m})\,d\beta^1_{t_{q+m}}\cdots g(t_q-t_{q+1})\,d\beta^1_{t_{q+1}}\\
&&=m\int_0^{t_q}u^{(m-1)}_{t_{q+1}}(t_q-t_{q+1})g(t_q-t_{q+1})\,d\beta_{t_{q+1}}^1.
\end{eqnarray*}
Then letting $N:=\sum_{j=1}^{k_i}s^1_j+s^2_j$, $N_m:=N-m$ and defining similar to above $\sum_{N,k_i,r,m}:=\sum_{m=0}^{\lfloor \frac{N}{2}\rfloor}\sum_{r\in Q_{N}^m}\sum_{r_{1}=1}^{k_i}\cdots \sum_{r_{N_m}=1}^{k_i}\frac{\prod_{j=1}^{N_m}s_j^1\max(1,s_j^2)}{2^m}$
gives 
\begin{eqnarray*}
&&\prod_{j=1}^{k_i}u_{t_{k_i}}^{(s_j)}(0)\\
&&\quad\quad=\sum_{N,k_i,r,m}\int_0^{t_{k_i}}\cdots\int_0^{t_{{k_i}+N_m-1}}g^{r_{N_m}}(t_{{k_i}+[N_m-1]}-t_{{k_i}+N_m})\,d\beta_{t_{{k_i}+N_m}}^{r_{N_m}}\\
&&\quad\quad\quad\quad\quad \cdots g^{r_{2}}(t_{{k_i}+[1]}-t_{{k_i}+2})\,d\beta_{t_{{k_i}+1}}^{r_{2}}g^{r_{1}}(t_{{k_i}}-t_{{k_i}+1})\,d\beta_{t_{{k_i}+1}}^{r_{1}},
\end{eqnarray*}
where $t_{{k_i}+[x]}\in \{t_{{k_i}},t_{{k_i}+1},...,t_{{k_i}+x}\}$, and $g^r(\cdot)=g^2(\cdot)$ if $r^2\neq 0$ and $g^r(\cdot)=g(\cdot)$ else. Similarly $\beta_{t}^{r}=t$ if $r^2=0$ and $\beta_{t}^{r}=\beta_{t}^{1}$ else. Without loss of generality we will again look only at one generic term, i.e. we set (by abuse of notation)
\begin{eqnarray*}
&&\prod_{j=1}^{k_i}u_{t_{k_i}}^{(s_j)}(0):=\int_0^{t_{k_i}}\cdots\int_0^{t_{{k_i}+N_m-1}}g^{r_{N_m}}(t_{{k_i}+[N_m-1]}-t_{{k_i}+N_m})\,d\beta_{t_{{k_i}+N_m}}^{r_{N_m}}\\
&&\quad\quad\quad\quad\quad \cdots g^{r_{2}}(t_{{k_i}+[1]}-t_{{k_i}+2})\,d\beta_{t_{{k_i}+1}}^{r_{2}}g^{r_{1}}(t_{{k_i}}-t_{{k_i}+1})\,d\beta_{t_{{k_i}+1}}^{r_{1}}.
\end{eqnarray*}
Due to $r\in Q_N^m$ we have $\sum_{i=1}^{N_m}1_{r_i\neq 0}=m$ and hence $m$ of $\beta^{r_{1}},...,\beta^{r_{N_m}}$ are Brownian motions and accordingly for $l=m-(k-2)$ we have
\begin{eqnarray*}
&&\mathbb{E}[f_{0}(T)H_l(T)]=\int_0^{T}\cdots\int_0^{t_{{k_i}+N_m-1}}g^{r_{N_m}}(t_{{k_i}+[N_m-1]}-t_{{k_i}+N_m})\,dt_{{k_i}+N_m}\\
&&\cdots g^{r_{2}}(t_{{k_i}+[1]}-t_{{k_i}+2})\,dt_{{k_i}+1}g^{r_{1}}(t_{{k_i}}-t_{{k_i}+1})\,dt_{{k_i}+1}d_{t_{k_i}}\cdots\,d_{t_{1}}+...,
\end{eqnarray*}
with $m$ of $g^{r_{1}}(\cdot),...,g^{r_{N_m}}(\cdot)$ corresponding to $g^2(\cdot)$ and the remaining to $g(\cdot)$ and thus $\mathbb{E}[f_{0}(T)H_l(T)]$ is well defined. By letting 
\begin{equation*}
G_{N_m}(t_{k_i},...,t_{k_i+N_m-1};g):=\int_0^{t_{k_i+N_m-1}}g^{r_{N_m}}(t_{k_i+[N_m-1]}-t_{k_i+N_m})\,dt_{k_i+N_m},
\end{equation*}
and for $l=1,...,N_m-1$
\begin{equation*}
G_l(t_{k_i},...,t_{k_i+l-1};g):=\int_{0}^{t_{k_i+l-1}}G_{l+1}(t_{k_i},...,t_{k_i+l};g)g^{r_l}(t_{k_i+[l-1]}+t_{k_i+l})\,d_{t_{k_i+l}},
\end{equation*}
the claim thus follows if for $t_{k_i}\in [0,T]$, $\vert G_1(t_{k_i};g)-G_1(t_{k_i};g_n)\vert \rightarrow 0$ for $n\rightarrow \infty$ holds true. As $(t_{k_i},...,t_{k_i+N_m-1})\mapsto G_{N_m}(t_{k_i},...,t_{k_i+N_m-1};h)$ is (absolutely) continuous and hence bounded on $[0,T]^{N_m}$ whenever $h\in L^2([0,T])$ the claim follows if $\vert \int_0^Tg^{r_1}(T-t_{1})-g_n^{r_1}(T-t_{1})\,dt_{1}\vert$ vanishes as $n\rightarrow \infty$, which holds true since $r^1 \in \{1,2\}$. 
\end{proof}
\end{lemma}
We denote by $\pi_i^n$ the $i$-th Malliavin weight that corresponds to the expansion of $X_T^{\epsilon,n}$. Further we recall that the Malliavin weights are also well-defined for $X_T^{\epsilon}$ as all quantities are evaluated in $\epsilon=0$. 
\begin{proposition}\label{roughWTX}Let $f$ be either in $C^\infty_p$ or bounded. If $f$ is also Lipschitz continuous and $g\in L^2([0,T])$ then \eqref{logPrice} admits a vol-of-vol expansion of arbitrary order $p\geq 2$, i.e. 
(recall \eqref{WTX}) 
\begin{equation*}\label{WTXp}
\vert \mathbb{E}[f(X_T^\epsilon)]-\sum_{i=0}^p\frac{\epsilon^i}{i!}\mathbb{E}[f(X_T^0)\pi_i]\vert=o(\epsilon^p),\quad \text{as } \epsilon \rightarrow 0.
\end{equation*}
\begin{proof}It follows from Lemma \ref{hdens} that there is a sequence $(g_n)\subset H$ such that $g_n\rightarrow g$ in $L^2([0,T])$. As $(g_n) \subset H$ we know that $X_T^{\epsilon, n}$ satisfies 
\begin{equation*}
\vert \mathbb{E}[f(X_T^{\epsilon,n})]-\sum_{i=0}^p\frac{\epsilon^i}{i!}\mathbb{E}[f(X_T^0)\pi^n_i]\vert=o(\epsilon^p),\quad \text{as } \epsilon \rightarrow 0.
\end{equation*}
Hence the claim follows from the Moore-Osgood Theorem if 
\begin{equation*}
\Big\vert \mathbb{E}[f(X_T^{\epsilon,n})]-\mathbb{E}[f(X_T^\epsilon)]+ \sum_{i=0}^p\frac{\epsilon^i}{i!}\Big(\mathbb{E}[f(X_T^0)\pi^n_i]-\mathbb{E}[f(X_T^0)\pi_i]\Big)\Big\vert=o(1),
\end{equation*}
as $n \rightarrow \infty$ uniformly in $\epsilon \in [-1,1]$. From Lemma \ref{Xg} it suffices to show 
\begin{equation*}
\vert \mathbb{E}[f(X_T^0)\pi^n_i]-\mathbb{E}[f(X_T^0)\pi_i] \vert=o(1) \quad \text{ as } n \rightarrow \infty
\end{equation*}
for $i=0,...,p$, but this follows from Lemma \ref{GGN} as from Proposition \ref{expREep}
\begin{eqnarray*}
&&\vert \mathbb{E}[f(X_T^0)\pi^n_i]-\mathbb{E}[f(X_T^0)\pi_i] \vert\\
&&\quad\quad\quad\leq \sum_{k=1}^p\sum_{l=k}^{p+2k}\vert \frac{\partial^l}{\partial x^l}\mathbb{E}[f(X_T^0)] \vert\vert \mathbb{G}_{k,l}(g_n)(T;i)-\mathbb{G}_{k,l}(g)(T;i)\vert,
\end{eqnarray*}
vanishes as $n\rightarrow \infty$.
\end{proof}
\end{proposition}

\bibliography{mybib2}{}

\providecommand{\bysame}{\leavevmode\hbox to3em{\hrulefill}\thinspace}
\providecommand{\MR}{\relax\ifhmode\unskip\space\fi MR }
\providecommand{\MRhref}[2]{%
  \href{http://www.ams.org/mathscinet-getitem?mr=#1}{#2}
}
\providecommand{\href}[2]{#2}
\begin{thebibliography}{10}

\bibitem{abi2019lifting}
Eduardo Abi~Jaber, \emph{Lifting the heston model}, Quantitative Finance
  (2019), 1--19.

\bibitem{OA2016}
Ozan Akdo\u{g}an, \emph{Variance curve models - finite dimensional realizations
  and beyond}, PhD Thesis, ETH Z\"{u}rich, 2016.

\bibitem{alos2012decomposition}
Elisa Al{\`o}s, \emph{A decomposition formula for option prices in the heston
  model and applications to option pricing approximation}, Finance and
  Stochastics \textbf{16} (2012), no.~3, 403--422.

\bibitem{alos2018exponentiation}
Elisa Alos, Jim Gatheral, and Rado{\v{s}} Radoi{\v{c}}i{\'c},
  \emph{Exponentiation of conditional expectations under stochastic
  volatility}, Quantitative Finance (2019), 1--15.

\bibitem{alos2007short}
Elisa Al{\`o}s, Jorge~A Le{\'o}n, and Josep Vives, \emph{On the short-time
  behavior of the implied volatility for jump-diffusion models with stochastic
  volatility}, Finance and Stochastics \textbf{11} (2007), no.~4, 571--589.

\bibitem{bayer2016pricing}
Christian Bayer, Peter Friz, and Jim Gatheral, \emph{Pricing under rough
  volatility}, Quantitative Finance \textbf{16} (2016), no.~6, 887--904.

\bibitem{bergomi2012stochastic}
Lorenzo Bergomi and Julien Guyon, \emph{Stochastic volatility's orderly
  smiles}, Risk \textbf{25} (2012), no.~5, 60.

\bibitem{ct2007}
Ren{\'e} Carmona and Michael Tehranchi, \emph{Interest rate models: an infinite
  dimensional stochastic analysis perspective}, Springer Science \& Business
  Media, 2007.

\bibitem{da2014stochastic}
Giuseppe Da~Prato and Jerzy Zabczyk, \emph{Stochastic equations in infinite
  dimensions}, Cambridge university press, 2014.

\bibitem{DF2001}
Damir Filipovi{\'c}, \emph{Consistency problems for heath-jarrow-morton
  interest rate models}, Springer, Lecture Notes in Mathematics, 1760, 2001.

\bibitem{friz2015large}
Peter~K Friz, Jim Gatheral, Archil Gulisashvili, Antoine Jacquier, and Josef
  Teichmann, \emph{Large deviations and asymptotic methods in finance}, vol.
  110, Springer, 2015.

\bibitem{fukasawa2011asymptotic}
Masaaki Fukasawa, \emph{Asymptotic analysis for stochastic volatility:
  martingale expansion}, Finance and Stochastics \textbf{15} (2011), no.~4,
  635--654.

\bibitem{gatheral2018volatility}
Jim Gatheral, Thibault Jaisson, and Mathieu Rosenbaum, \emph{Volatility is
  rough}, Quantitative Finance \textbf{18} (2018), no.~6, 933--949.

\bibitem{gatheral2018affine}
Jim Gatheral and Martin Keller-Ressel, \emph{Affine forward variance models},
  Finance and Stochastics (2019), 1--33.

\bibitem{gerasimovics2018h}
Andris Gerasimovics and Martin Hairer, \emph{H{\"o}rmander's theorem for
  semilinear spdes}, arXiv preprint arXiv:1811.06339 (2018).

\bibitem{gulisashvili2019higher}
Archil Gulisashvili, Ra{\'u}l Merino, Marc Lagunas, and Josep Vives,
  \emph{Higher order approximation of call option prices under stochastic
  volatility models}, arXiv preprint (2019).

\bibitem{harms2019strong}
Philipp {Harms}, \emph{{Strong convergence rates for Markovian representations
  of fractional Brownian motion}}, arXiv e-prints (2019).

\bibitem{KT2003}
Naoto Kunitomo and Akihiko Takahashi, \emph{On validity of the asymptotic
  expansion approach in contingent claim analysis}, The Annals of Applied
  Probability \textbf{13} (2003), no.~3, 914--952.

\bibitem{lewis2000option}
Alan Lewis, \emph{Option valuation under stochastic volatility}, Tech. report,
  Finance Press, 2000.

\bibitem{musiela1993stochastic}
Marek Musiela, \emph{Stochastic pdes and term structure models}, Journ{\'e}es
  Internationales de Finance (1993).

\bibitem{nd2006}
David Nualart, \emph{The malliavin calculus and related topics}, Springer,
  2006.

\bibitem{PP2005}
Philip~E Protter, \emph{Stochastic differential equations}, Stochastic
  Integration and Differential Equations, Springer, 2005, pp.~249--361.

\bibitem{TS2011}
Maria Siopacha and Josef Teichmann, \emph{Weak and strong taylor methods for
  numerical solutions of stochastic differential equations}, Quantitative
  Finance \textbf{11} (2011), no.~4, 517--528.

\bibitem{ttt2012}
Akihiko Takahashi, Kohta Takehara, and Masashi Toda, \emph{A general
  computation scheme for a high-order asymptotic expansion method},
  International Journal of Theoretical and Applied Finance \textbf{15} (2012),
  no.~06, 1250044.

\bibitem{TT2012}
Akihiko Takahashi and Toshihiro Yamada, \emph{An asymptotic expansion with
  push-down of malliavin weights}, SIAM Journal on Financial Mathematics
  \textbf{3} (2012), no.~1, 95--136.

\bibitem{ty2004}
Akihiko Takahashi and Nakahiro Yoshida, \emph{An asymptotic expansion scheme
  for optimal investment problems}, Statistical Inference for Stochastic
  Processes \textbf{7} (2004), no.~2, 153--188.

\bibitem{WS1987}
Shinzo Watanabe, \emph{Analysis of wiener functionals (malliavin calculus) and
  its applications to heat kernels}, The annals of Probability (1987), 1--39.

\bibitem{YN1992}
Nakahiro Yoshida, \emph{Asymptotic expansions of maximum likelihood estimators
  for small diffusions via the theory of malliavin-watanabe}, Probability
  Theory and Related Fields \textbf{92} (1992), no.~3, 275--311.

\end{thebibliography}
\bibliographystyle{amsplain}

\end{document}